\documentclass[11pt]{article}

\makeatletter
\def\@begintheorem#1#2{\trivlist%
 \item[\hskip \labelsep{\sffamily\bfseries #2\ #1}]\itshape}
\newtheorem{teo}{Theorem}[section]
\newtheorem{defi}[teo]{Definition}
\newtheorem{cor}[teo]{Corollary}
\newtheorem{lem}[teo]{Lemma}
\newtheorem{pro}[teo]{Proposition}
\newtheorem{_rem}[teo]{Remark}
\newtheorem{_eje}[teo]{Example}
\newenvironment{rem}{\def\@begintheorem##1##2{\trivlist%
 \item[\hskip\labelsep{\sffamily\bfseries ##2\ ##1}]}\begin{_rem}}{\end{_rem}}

\makeatother

\newenvironment{beweis}{{\em Proof:}}{\hfill $\rule{2mm}{2mm}$
\vspace{3mm}

}

\DeclareMathAlphabet{\Ma}{U}{msa}{m}{n}
\DeclareMathAlphabet{\Mb}{U}{msb}{m}{n}
\DeclareMathAlphabet{\Meuf}{U}{euf}{m}{n}

\def\got#1{\Meuf{#1}}

\def\mr #1.{\mathrm{#1}}

\DeclareSymbolFont{ASMa}{U}{msa}{m}{n}
\DeclareSymbolFont{ASMb}{U}{msb}{m}{n}
\DeclareMathSymbol{\hrist}{\mathord}{ASMa}{"16}
\DeclareMathSymbol{\varkappa}{\mathalpha}{ASMb}{"7B}
\DeclareMathSymbol{\CrPr}{\mathord}{ASMb}{"6F}

\newfont{\EinsFont}{cmr7 scaled 1070}
\def\EINS{{\mathchoice{
 \mbox{\unitlength1cm\begin{picture}(.25,.2)\put(0,0){$1$}%
 \put(0.105,0){{\mbox{\fontfamily{cmr}\upshape\small l}}}\end{picture}}}{%
 \mbox{\unitlength1cm\begin{picture}(.25,.2)\put(0,0){$1$}%
 \put(0.105,0){{\mbox{\fontfamily{cmr}\upshape\small l}}}\end{picture}}}{%
 \mbox{\unitlength1cm\begin{picture}(.18,.15)\put(0,0){$\scriptstyle 1$}%
 \put(0.07,0){{\mbox{\fontfamily{cmr}\upshape\EinsFont l}}}\end{picture}}}{%
 \mbox{\unitlength1cm\begin{picture}(.18,.15)\put(0,0){$\scriptstyle 1$}%
 \put(0.07,0){{\mbox{\fontfamily{cmr}\upshape\EinsFont l}}}\end{picture}}}}}

\def\restriction{{\mathchoice{
 \mbox{\unitlength1cm\begin{picture}(.2,.4)%
  \bezier{5}(.07,.3)(.1,.27)(.13,.24)%
  \put(.07,.35){\line(0,-1){.5}}\end{picture}}}{
 \mbox{\unitlength1cm\begin{picture}(.2,.4)%
  \bezier{5}(.07,.3)(.1,.27)(.13,.24)%
  \put(.07,.35){\line(0,-1){.5}}\end{picture}}}{
  \hrist}{\hrist}}}

  \def\al #1.{{\mathcal{#1}}}
  \def\ot #1.{{\got{#1}}}
  \def\C{\Mb{C}}
  
\def\ed{\end{document}}
\def\be{\begin{equation}}
\def\ee{\end{equation}}
\def\bea{\begin{eqnarray}}
\def\eea{\end{eqnarray}}
\def\beaO{\begin{eqnarray*}}
\def\eeaO{\end{eqnarray*}}

\newcommand{\kSPgr}{\mbox{$\seins_{\kern-.2em\szwei}^{\kern-.1em\sdrei}$}}

\DeclareMathSymbol{\hsemi}{\mathord}{ASMb}{"6E}
\newcommand{\semi}[2]{\mbox{$#1\kern.1em\hsemi\kern.1em#2$}}

\def\LA{\left\langle\bgroup}
\def\LE{\left[\bgroup}
\def\LG{\left\{\bgroup}
\def\LR{\left(\bgroup}
\def\RA{\egroup^{\rule{0mm}{2mm}}\right\rangle}
\def\RE{\egroup^{\rule{0mm}{2mm}}\right]}
\def\RG{\egroup^{\rule{0mm}{2mm}}\right\}}
\def\RR{\egroup^{\rule{0mm}{2mm}}\right)}
\def\Ldummy{\left.\bgroup}
\def\Rdummy{\egroup^{\rule{0mm}{2mm}}\right.}

\def\Kbegin{\begin{equation} \left. \begin{array}{rcl}}
\def\Kend{\end{array} \right\} \end{equation}}

\def\l2{\Lambda^{\mbox{\tiny $(2)$}}}

\title{\bf An application of the DR-duality theory for compact
groups to endomorphism categories of C*-algebras with nontrivial
center}
\author{
 {\sc Hellmut Baumg\"artel}  \\[2mm] 
 {\footnotesize Mathematical Institute, University of Potsdam,}     \\   
 {\footnotesize Am Neuen Palais 10, Postfach 601~553,}             \\ 
 {\footnotesize D--14415 Potsdam, Germany.}                        \\
 {\footnotesize baumg@rz.uni-potsdam.de}
\and    
 {\sc Fernando Lled\'o }\thanks{Present address: {\em
  Institute of Pure and Applied Mathematics, RWTH--Aachen,
  Templergraben 55, D--52056 Aachen, Germany.} 
  E-Mail:{\em lledo@iram.rwth-aachen.de}}  \\[2mm] 
 {\footnotesize Max--Planck--Institut f\"ur Gravitationsphysik,}     \\   
 {\footnotesize Albert--Einstein--Institut,}
 {\footnotesize Am M\"uhlenberg 1,}             \\ 
 {\footnotesize D--14476 Golm, Germany.}                        \\
 {\footnotesize lledo@aei-potsdam.mpg.de}}

\date{\today{}}

\begin{document}
\maketitle
\begin{center}
{\sl Dedicated to Sergio Doplicher and John E.~Roberts  
 on the occasion of their 60th birthdays.}
\end{center}
\vspace{.5cm}

\begin{abstract}
The main result of the Doplicher/Roberts (DR--) duality theory for
compact groups, applied to the special case of endomorphism categories
with permutation and conjugation structure of a fixed C*--algebra
${\cal A}$ with trivial center, says that such a category can be
characterized as the category of all canonical endomorphisms of
${\cal A}$ w.r.t.~an (essentially uniquely determined) Hilbert
extension $\{{\cal F},{\cal G}\}$ of ${\cal A}$, where ${\cal
G}$ is a compact automorphism group of ${\cal F}$ and
${\cal A}'\cap {\cal F}= \Mb{C}\EINS.$

In \cite{Lledo97b} C*--Hilbert systems 
$\{{\cal F},{\cal G}\}$
are considered where the fixed point
algebra ${\cal A}$ has nontrivial center ${\cal Z}$ and where 
${\cal A}'\cap {\cal F}={\cal Z}$
is satisfied. 
The corresponding category of all canonical endomorphisms of
${\cal A}$
contains characteristic mutually isomorphic subcategories of the
DR--type which are connected with the choice of distinguished 
${\cal G}$--invariant algebraic Hilbert spaces within the corresponding
${\cal G}$--invariant Hilbert 
${\cal Z}$--modules.

We present in this paper the solution of the corresponding inverse
problem. More precisely, assuming that the given endomorphism
category ${\cal T}$
of a C*--algebra ${\cal A}$ with center $\al Z.$
contains a certain subcategory of the DR--type, a Hilbert extension
$\{{\cal F},{\cal G}\}$ of ${\cal A}$ is constructed such that 
${\cal T}$
is isomorphic to the category of all canonical endomorphisms of
${\cal A}$ w.r.t.~$\{{\cal F},{\cal G}\}$ and
${\cal A}'\cap {\cal F}={\cal Z}$.
Furthermore, there is a natural equivalence relation between
admissible subcategories and it is shown that two admissible
subcategories yield
${\cal A}$--module isomorphic Hilbert extensions iff they are equivalent.
The essential step of the solution is the application of the
standard DR--theory to the assigned subcategory.
\end{abstract}
\begin{center}
MSC--class: 47L65, 22D25, 46L08
\end{center}

\section{Introduction}
One of the origins of the DR--duality theory for compact groups
(cf.~\cite{Doplicher89b}) is the analysis of superselection 
structures formulated in the context of
Algebraic Quantum Field Theory, where endomorphism categories
${\cal T}$ of a fixed C*--algebra ${\cal A}$
with trivial center (the so--called algebra of 
quasilocal observables) appear. These categories ${\cal T}$ are
equipped with a conjugation and a permutation structure.
Application of the general DR--duality theory yields to the result
that ${\cal T}$ can be characterized by a compact symmetry group
${\cal G}$ via the construction of a suitable Hilbert extension
$\{{\cal F},{\cal G}\}$ of ${\cal A}$, where ${\cal F}$ is
called the field algebra 
(see \cite{Doplicher89a,Doplicher90,Baumgaertel97}). The condition
${\cal A}'\cap {\cal F}=\Mb{C}\EINS$ is crucial for this construction.

This result suggests the natural question for 
conditions of endomorphism categories ${\cal T}$
of C*--algebras ${\cal A}$ with nontrivial center 
${\cal Z}$ such that there is still a description of ${\cal T}$
by a compact group. Already the paper \cite[Sections~2,3]{Doplicher89b} 
starts with nontrivial center.
The point there is that the ``flip property" for the permutator
$\epsilon(\alpha,\beta)$
and the intertwiners
$A\in (\alpha,\alpha'),\,B\in (\beta,\beta')$, i.e. 
\[
\epsilon(\alpha',\beta')A\times B=B\times A\,\epsilon(\alpha,\beta),
\]
is assumed to be valid for {\em all} intertwiners. Only later on
the condition $(\iota,\iota)=\Mb{C}\EINS$
is added to arrive at the famous DR--theorem.
Also in the context of more general categories (that do not 
assume the existence of a permutator) some results for 
nontrivial $(\iota,\iota)$ are stated (cf.~\cite[Section~2]{Longo97}).

In \cite{Lledo97b} we present
several results on general Hilbert C*--systems, where
${\cal A}'\cap {\cal F}={\cal Z}\supset \Mb{C}\EINS$
is assumed. For example, properties of the corresponding
category of canonical endomorphisms are proved and the breakdown
of the Galois correspondence between the symmetry group
${\cal G}$
and the stabilizer of its fixed point algebra
${\cal A}$
is stated. Further a new concept of ``irreducibility" is proposed.
In \cite{BaumgaertelIn99} further properties of ``irreducible
endomorphisms" are mentioned and for special automorphism
categories the solution of the ``inverse problem" (given the
category, construct an assigned Hilbert extension) is briefly
described. 
The paper \cite{pLledo00} contains a description of the
status of the problem for automorphism categories
as a counterpart of the special case where 
${\cal Z}=\Mb{C}\EINS$. In the previous two references the 
corresponding group ${\cal G}$ is abelian.

The present paper describes the solution of the inverse
problem in the general case. First, sufficient conditions for the
endomorphism category ${\cal T}$ of ${\cal A}$
are stated such that there exists a Hilbert extension
$\{{\cal F},{\cal G}\}$
of ${\cal A}$ and ${\cal T}$ turns out to be the category
of all canonical endomorphisms of
${\cal A}$. The crucial assumption is the existence of a certain subcategory
${\cal T}_{\Mb{C}}$ of ${\cal T}$
of the DR--type. The goal of this note is to show that the
DR--properties of ${\cal T}_{\Mb{C}}$
are sufficient to characterize the category
${\cal T}$ as the category of all canonical endomorphisms of
${\cal A}$ w.r.t.~a suitable Hilbert extension
$\{{\cal F},{\cal G}\}$ of ${\cal A}$, where ${\cal G}$
is the characteristic compact DR--group for ${\cal T}_{\Mb{C}}$,
interpreted as an automorphism group of the Hilbert extension.
The conditions are also necessary. This is an easy implication
of the results for Hilbert C*--systems, obtained in \cite{Lledo97b}. 
Second, a uniqueness result is stated: Two admissible subcategories of
${\cal T}$ yield ${\cal A}$--module isomorphic Hilbert extensions iff 
they are equivalent (in a precise sense formulated in the following
section).

The present paper is structured in 4 sections: In the following section
we present the postulates for the category $\al T.$ that assure the 
unique (up to $\al A.$--module isomorphy) extension of the 
C*--algebra $\al A.$ (cf.~Theorem~\ref{Main}, which is our main result). 
In Subsection~\ref{S3.1} we introduce the new notion of irreducible 
endomorphism in $\al T.$ and show some of its consequences. Finally,
we give in the following two subsections the proof of 
Theorem~\ref{Main}. 

For convience of the reader we finish this section recalling
some useful definitions concerning Hilbert extensions
of a C*--algebra $\al A.$: 
A C*--system $\{{\cal F},{\cal G}\}$, where 
${\cal G}\subset\mathrm{aut}\,{\cal F}$
is compact w.r.t.~the topology of pointwise norm convergence, is called
{\em Hilbert} if $\mathrm{spec}\,{\cal G}=\hat{\cal G}$ 
(the dual of ${\cal G}$) and if each spectral subspace 
${\cal F}_{D},\,D\in\hat{\cal G}$, contains an algebraic Hilbert space
${\cal H}_{D}$ such that ${\cal G}$ acts invariantly on ${\cal H}_{D}$
and the unitary representation ${\cal G}\restriction {\cal H}_{D}$
is a element of the equivalence class $D$.
A Hilbert space ${\cal H}\subset {\cal F}$ is called {\em algebraic}
if the scalar product $\langle\cdot,\cdot\rangle$ of
${\cal H}$ is given by $\langle A,B\rangle\EINS:=A^{\ast}B,\,A,B\in{\cal H}.$
A Hilbert system $\{{\cal F},{\cal G}\}$ is called a  {\em Hilbert extension}
of a C*--algebra ${\cal A}$ if ${\cal A}\subset {\cal F}$ is the 
fixed point algebra of ${\cal G}$. Two Hilbert extensions 
$\{{\cal F}_1,{\cal G}\}$,$\{{\cal F}_2,{\cal G}\}$ are called
$\al A.$--module isomorphic if the exists an isomorphism 
$\al J.\colon\al F._1\to\al F._2$, such that $\al J.(A)=A$, $A\in\al A.$,
that intertwines between the corresponding group actions.
Finally, the {\em stabilizer} $\mathrm{stab}\,{\cal A}$
for a unital C*--subalgebra ${\cal A}\subset{\cal F}$ is defined by
$\mathrm{stab}\,{\cal A}:=\{g\in\mathrm{aut}\,{\cal F}\mid g(A)=A
\,\mbox{for all}\,A\in{\cal A}\}.$
(The terminology Hilbert system can be traced back to \cite{Doplicher88}, 
where in the case just mentioned the spectrum of ${\cal G}$
is called the Hilbert spectrum.)

\section{Assumptions on the endomorphism category and the main result}
\label{Assu}

In the present section we will collect the assumptions on the category
${\cal T}$ of suitable endomorphisms of a unital C*--algebra
${\cal A}$ that guarantees the main theorem stated at the end of this
section. For standard notions within category theory we refer to
\cite{bMacLane98}.

Let ${\cal T}$ be a tensor C*--category of unital endomorphisms $\alpha$
of ${\cal A}$ \cite{Doplicher89b}.
We denote the objects by
$\alpha,\beta,\gamma,\ldots\in \mathrm{Ob}\,{\cal T}$.
The arrows between objects $\alpha,\beta$
are given as usual by the intertwiner spaces
$(\alpha,\beta):=\{X\in {\cal A}\mid X\alpha(A)=\beta(A)X,\,A\in
{\cal A}\}$ and we put $A\times B:= A\alpha(B)$
for $A\in (\alpha,\alpha'), B\in (\beta,\beta')$, so 
that $A\times B\in (\alpha\beta,\alpha'\beta')$.
By $\iota$ we denote the identical endomorphism and
$(\iota,\iota)={\cal Z}$ is the center of ${\cal A}$.
Note that $(\alpha,\beta)$ is a left 
$\beta({\cal Z})$-- and a right $\alpha({\cal Z})$--module, 
i.e.~$\beta({\cal Z})(\alpha,\beta)\alpha({\cal
Z})=(\alpha,\beta)\alpha({\cal Z})\subseteq (\alpha,\beta)$.
The conditions on $\al T.$ are given by:
 
\begin{itemize}
\item[P.1.1] 
  ${\cal T}$ is closed w.r.t.~direct sums
  $\alpha\oplus\beta$, i.e.~if
  $\alpha,\beta\in\mathrm{Ob}\,{\cal T}$, then there are isometries
  $V,W\in{\cal A}$ with
  $V^{\ast}W=0,\,VV^{\ast}+WW^{\ast}=\EINS$ such that
  $\gamma(\cdot):=
  V\alpha(\cdot)V^{\ast}+W\beta(\cdot)W^{\ast}\in\mathrm{Ob}\,{\cal T}.$
  In this case $V\in (\alpha,\gamma),W\in (\beta,\gamma)$.
\item[P.1.2] 
  ${\cal T}$ is closed w.r.t.~subobjects $\beta < \alpha$,
  i.e.~if $\alpha\in\mathrm{Ob}\,{\cal T}$
  and $\beta$ is a unital endomorphism of ${\cal A}$
  such that there is an isometry $V\in (\beta,\alpha)$,
  then $\beta\in\mathrm{Ob}\,{\cal T}$.
  In this case $\beta(\cdot)=V^{\ast}\alpha(\cdot)V.$

\item[P.1.3]
  ${\cal T}$ is closed w.r.t.~complementary subobjects, i.e.~if
  $\alpha\in\mathrm{Ob}\,{\cal T}$ and $\beta <\alpha$,
  then there is a subobject $\beta' < \alpha$
  such that $\alpha=\beta\oplus\beta'$.

\item[P.2]
  ${\cal T}$ contains a C*--subcategory ${\cal T}_{\Mb{C}}$
  with $\mathrm{Ob}\,{\cal T}_{\Mb{C}}=\mathrm{Ob}\,{\cal T},$
  where the arrows $(\alpha,\beta)_{\Mb{C}}\subset (\alpha,\beta)$
  satisfy the following properties:
\begin{itemize}
\item[P.2.1] $(\beta,\gamma)_{\Mb{C}}\cdot
            (\alpha,\beta)_{\Mb{C}}\subseteq (\alpha,\gamma)_{\Mb{C}},$
\item[P.2.2] $\alpha(\beta,\gamma)_{\Mb{C}}\subseteq(\alpha\beta,
            \alpha\gamma)_{\Mb{C}},$
\item[P.2.3] $(\alpha,\beta)_{\Mb{C}}
            \subseteq(\alpha\gamma,\beta\gamma)_{\Mb{C}}$,
\item[P.2.4] 
            $(\alpha,\beta)_{\Mb{C}}^{\ast}\subseteq (\beta,\alpha)_{\Mb{C}}$,
\item[P.2.5] $(\iota,\iota)_{\Mb{C}}=\Mb{C}\EINS,$
\item[P.2.6] Any finite set of linearly independent elements 
            $F_{1},F_{2},\ldots,F_{n}\in (\alpha,\beta)_{\Mb{C}}$
            (a complex Banach space), is linearly independent 
            modulo $\alpha({\cal Z})$ in $(\alpha,\beta)$, i.e.~if
            $\sum_{j=1}^{n}\lambda_{j}F_{j}=0,\,\lambda_{j}\in\Mb{C}$,
            implies $\lambda_{j}=0$,
            then also $\sum_{j=1}^{n}F_{j}\alpha(Z_{j})=0$
            implies $\alpha(Z_{j})=0$.
            Moreover, $(\alpha,\beta)_{\Mb{C}}$
            generates $(\alpha,\beta)$, i.e.
            $(\alpha,\beta)=(\alpha,\beta)_{\Mb{C}}\alpha({\cal Z})=
            \beta({\cal Z})(\alpha,\beta)_{\Mb{C}}\alpha({\cal Z})$.
\item[P.2.7] ${\cal T}_{\Mb{C}}$ is closed w.r.t.~direct sums, 
            subobjects and complementary subobjects, i.e.~now the
            required projections and isometries must ly in the corresponding
            intertwiner spaces $(\cdot,\cdot)_{\Mb{C}}$.
\end{itemize}
\item[P.3] There is a permutation structure\footnote{In 
  \cite{Doplicher89b} categories
  as ${\cal T}_{\Mb{C}}$ satisfying the properties below are
  called {\em symmetric}.}, on ${\cal T}_{\Mb{C}}$,
  i.e.~there is a mapping
  $\{\alpha,\beta\}\rightarrow \epsilon(\alpha,\beta)\in
  (\alpha\beta,\beta\alpha)_{\Mb{C}}$, where $\epsilon(\alpha,\beta)$
  is a unitary satisfying:
\begin{itemize}
\item[P.3.1] $\epsilon(\alpha,\beta)\cdot\epsilon(\beta,\alpha)=\EINS,$
\item[P.3.2] $\epsilon(\iota,\alpha)=\epsilon(\alpha,\iota)=\EINS,$
\item[P.3.3] $\epsilon(\alpha\beta,\gamma)=\epsilon(\alpha,\gamma)\cdot
            \alpha(\epsilon(\beta,\gamma)),$
\item[P.3.4] $\epsilon(\alpha',\beta')A\times B=
            B\times A\,\epsilon(\alpha,\beta)\,
            \mbox{for all}\,A\in(\alpha,\alpha')_{\Mb{C}},
            B\in(\beta,\beta')_{\Mb{C}}.$
\end{itemize}

\item[P.4] There is a conjugation structure on ${\cal T}_{\Mb{C}}$,
i.e.~to each $\alpha\in\mathrm{Ob}\,{\cal T}_{\Mb{C}}$
there corresponds a conjugate
$\overline{\alpha}\in\mathrm{Ob}\,{\cal T}$ and intertwiners
$R_{\alpha}\in(\iota,\overline{\alpha}\alpha)_{\Mb{C}},\,
S_{\alpha}\in(\iota,\alpha\overline{\alpha})_{\Mb{C}}$ such that
\begin{itemize}
\item[P.4.1] $S_{\alpha}^{\ast}\alpha(R_{\alpha})=\EINS,\quad 
            R_{\alpha}^{\ast}\overline{\alpha}(S_{\alpha})=\EINS,$
\item[P.4.2] $S_{\alpha}=\epsilon(\overline{\alpha},\alpha)R_{\alpha}.$
\end{itemize}
\end{itemize}

\begin{rem}\label{Comm}
\begin{itemize}
\item[(i)]
  The preceding axioms imply that the subcategory ${\cal T}_{\Mb{C}}$
  satisfies the postulates of the DR--theory developed in \cite{Doplicher89b}, 
  so that we can apply standard results from this theory. 
\item[(ii)]  Note that the definition of subobject in (P.1.2) is 
  {\em not} a straightforward generalization to nontrivial center 
  situation of the one given in \cite{Doplicher89b}. Namely,
  if $E$ is a selfadjoint central projection $E\in\al Z.$ with
  $0<E<\EINS$, then there is no isometry $W\in\al A.$ such that 
  $E=WW^*$, because in that case we simply have
\[
 \EINS=W^*W\,W^*W=W^*\,E\,W=E\,W^*W=E\,.
\]
\item[(iii)] From (P.4) it already follows that
  ${\cal T}$ has conjugates (cf.~\cite[Section~2]{Longo97}).
  Further (P.2.7) also shows that $\al T.$ is closed under direct sums,
  subobjects and complementary subobjects.
  (P.2.1)--(P.2.3) imply
  $(\alpha,\alpha')_{\Mb{C}} \times (\beta,\beta')_{\Mb{C}}\subseteq 
  (\alpha\beta, \alpha'\beta')_{\Mb{C}}$ and
  from (P.2.4) it follows immediately that the equality
  $(\alpha,\beta)_{\Mb{C}}^{\ast}=(\beta,\alpha)_{\Mb{C}}$ holds.
\item[(iv)] (P.1.1) implies an ``a priori property" of ${\cal A}$,
  namely there are two isometries
  $V,W\in {\cal A}, V^{\ast}W=0, VV^{\ast}+WW^{\ast}=\EINS$.
  (P.1.3) implies that if $E\in(\alpha,\alpha)$
  is a projection such that there is an isometry $V$
  with $VV^{\ast}=E$, then there is also an isometry
  $W$ with $WW^{\ast}=\EINS-E.$
\end{itemize}   
\end{rem}

\begin{rem}
Let $\mathrm{Ob}\,{\cal T}\ni\alpha\rightarrow 
V_{\alpha}\in (\alpha,\alpha)$ 
be a choice of unitaries that satisfy
\begin{equation}\label{Cond}
V_{\alpha\circ\beta}=V_{\alpha}\times V_{\beta}.
\end{equation}
Note that (\ref{Cond}) implies $V_{\iota}=\EINS$, because 
$V_{\beta}=V_{\iota\beta}=V_{\iota}\times V_{\beta}=V_{\iota}V_{\beta}$.
This choice allows to define from
the subcategory ${\cal T}_{\Mb{C}}$ of ${\cal T}$
another subcategory ${\cal T}_{\Mb{C}}'$ of ${\cal T}$
satisfying the same properies as ${\cal T}_{\Mb{C}}$.
Indeed, put
\begin{equation}\label{Int}
(\alpha,\beta)_{\Mb{C}}':=V_{\beta}(\alpha,\beta)_{\Mb{C}}V_{\alpha}^{\ast}
 \subset (\alpha,\beta)
\end{equation}
and the corresponding permutation structure $\epsilon'(\cdot,\cdot)$
for ${\cal T}_{\Mb{C}}'$ is given by
\begin{equation}\label{Ep}
\epsilon'(\alpha,\beta):=V_{\beta}\times V_{\alpha}\cdot
\epsilon(\alpha,\beta)\cdot (V_{\alpha}\times V_{\beta})^{\ast}.
\end{equation}
It is easy to check that 
$\epsilon'$ satisfies the properties in (P.3). The corresponding conjugates
$R_{\alpha}'$ are defined by
\begin{equation}\label{RS}
R_{\alpha}':=V_{\overline{\alpha}\alpha}R_{\alpha}\quad\mathrm{and}\quad
S_{\alpha}':=\epsilon'(\alpha,\overline{\alpha})R_{\alpha}'.
\end{equation}
Then it is straightforward to verify the postulates (P.2)--(P.4)
for the new subcategory.
\end{rem}

The preceding remark suggests to define an equivalence relation between
different subcategories of ${\cal T}$:

\begin{defi}
The subcategories ${\cal T}_{\Mb{C}}$
and ${\cal T}_{\Mb{C}}'$,
satisfying the postulates (P.2)--(P.4),
are called equivalent, if there is an assignment
\[
\mathrm{Ob}\,{\cal T}\ni\alpha\rightarrow V_{\alpha}\in(\alpha,\alpha),\;\,
\mbox{with}\; V_{\alpha}\;\mbox{unitary~and}\;\;
V_{\alpha\circ\beta}=V_{\alpha}\times V_{\beta}\,,
\]
such that Eqs.~(\ref{Int}),(\ref{Ep}) and (\ref{RS}) hold.
\end{defi}

\begin{teo}\label{Main}
\begin{itemize}
\item[(i)] Let ${\cal T}$  satisfy the postulates (P.1)--(P.4) before. 
      Then there exists a Hilbert extension $\{{\cal F},{\cal G}\}$
      of ${\cal A}$ with ${\cal A}'\cap{\cal F}={\cal Z}$
      such that ${\cal T}$  is isomorphic to the category of all canonical
      endomorphisms of $\{{\cal F},{\cal G}\}$.
\item[(ii)] Further let ${\cal T}_{\Mb{C}},{\cal T}_{\Mb{C}}'$ 
      be two subcategories of ${\cal T}$ satisfying the postulates 
      (P.2)--(P.4). Then the corresponding Hilbert extensions are
      ${\cal A}$--module isomorphic iff
      ${\cal T}_{\Mb{C}}$ and ${\cal T}_{\Mb{C}}'$ are equivalent.
\end{itemize}
\end{teo}

\section{Proof of Theorem~\ref{Main}}

In the present section we will give a (constructive) proof of the
previous theorem. For this purpose we will use well--known
results already stated in 
\cite{Doplicher89b,Baumgaertel97,Lledo97b}
as well as in \cite[Sections~2,3]{Longo97}.

First we study properties of
${\cal T}$ implied by the existence of the subcategory 
${\cal T}_{\Mb{C}}$, in particular we introduce the notion of 
irreducibility in the context of ${\cal T}$
and prove the decomposition theorem.

\subsection{Irreducibility and decomposition theorem}\label{S3.1}

Since the C*--algebra ${\cal A}$
has a nontrivial center ${\cal Z}$
it is immediate that one needs to extend the notion of
irreducible objects to the category
${\cal T}$ (cf.~e.g.~\cite[Section 5]{Lledo97b}). We propose

\begin{defi}\label{3.1}
$\rho\in\mathrm{Ob}\,{\cal T}$ is called irreducible if
$(\rho,\rho)={\cal Z}$. 
We denote the set of all irreducible objects of
${\cal T}$ by $\mathrm{Irr}\,{\cal T}$ and by
$\mathrm{Irr}_{0}\,{\cal T}$ a complete system of irreducible and 
mutually disjoint objects of ${\cal T}$.
\end{defi}

Note that irreducibility of $\rho$
in the sense of Definition~\ref{3.1} and of
$\rho$ in the usual sense as an object in ${\cal T}_{\Mb{C}}$
coincide, because
$\rho\in\mathrm{Irr}\,{\cal T}$ iff
$(\rho,\rho)_{\Mb{C}}=\Mb{C}\EINS$. We state some further
consequences of this definition.

\begin{lem}
\begin{itemize} 
\item[(I)] If $\rho,\sigma\in\mathrm{Irr}\,{\cal T}$, then either
           $\rho$ is unitarily equivalent to $\sigma$ 
           or they are disjoint (i.e.~$(\rho,\sigma)=\{0\})$.

\item[(II)]The following properties are equivalent:
\begin{itemize}
\item[(i)] $\rho$ is irreducible,
\item[(ii)] $\overline{\rho}$ is irreducible,
\item[(iii)] $(\rho,\rho)=\rho({\cal Z})$,
\item[(iv)] $(\iota,\overline{\rho}\rho)=R_{\rho}{\cal Z}$, where
            $R_{\rho}$ is a conjugate according to property (P.4).
\end{itemize}
\item[(III)] If $\rho$ is irreducible, then
$(\rho,\alpha)_{\Mb{C}}$ is an algebraic Hilbert space in ${\cal A}$
for each $\alpha\in\mathrm{Ob}\,{\cal T}$ and
$(\rho,\alpha)=(\rho,\alpha)_{\Mb{C}}{\cal Z}$
is a right--${\cal Z}$--Hilbert module with the scalar product
$\langle X,Y\rangle :=X^{\ast}Y.$
\item[(IV)] $\rho$ is irreducible iff there is no subobject of $\rho$.
\end{itemize}
\end{lem}
\begin{beweis}
(I) Let $(\rho,\sigma)\supset \{0\}$, so that the inclusion
$(\rho,\sigma)_{\Mb{C}}\supset\{0\}$ is also proper.
Then it is straightforward to construct a unitary
$U\in (\rho,\sigma)_{\Mb{C}}$.

(II) The proof  uses the vector space isomorphisms between
intertwiner spaces (see, for example \cite[Lemma~2.1]{Longo97}) together
with the link between conjugates and permutation given by
assumption (P.4.2). Namely, the latter implies that
$(\alpha,\alpha)=\alpha({\cal Z})$ iff
$(\overline{\alpha},\overline{\alpha})=\overline{\alpha}({\cal Z})$
for all
$\alpha\in\mathrm{Ob}\,{\cal T}$, while the vector space 
isomorphisms yield $(\alpha,\alpha)={\cal Z}$ iff
$(\overline{\alpha},\overline{\alpha})=\overline{\alpha}({\cal Z}).$

(III) It follows immediately from (P.2.3)--(P.2.5).

(IV) If ${\rho}$ is irreducible, then it is straightforward 
to see that any isometry $W\in (\sigma,\rho)$ is actually
a unitary, because $WW^*=:E\in (\rho,\rho)$, hence 
$E=\EINS$ by Remark~\ref{Comm}~(ii).
To show the reverse implication assume that $\rho$
is not irreducible. Now the estimate
\[
A\leq (R_{\rho}^{\ast}R_{\rho})^{2}
\Phi_{\rho}(A),\quad 0\leq A \in (\rho,\rho)_{\Mb{C}}\,,
\]
where $\Phi_{\rho}$ denotes the corresponding left inverse
(see for example \cite[Lemma~2.7]{Longo97}), implies that the C*--algebra
$(\rho,\rho)_{\Mb{C}}\supset \Mb{C}\EINS$
is finite--dimensional, hence $\rho$
must have proper subobjects (see 
e.g.~\cite[Lemma~11.1.27 and 11.1.29]{bBaumgaertel92}).
The crucial fact is that
$\Phi_{\alpha}(A)\in\Mb{C}\EINS$
for $A\in (\alpha,\alpha)_{\Mb{C}}$.
\end{beweis}

The previous lemma implies in particular that the restriction of 
irreducible objects to ${\cal Z}$ are automorphisms of ${\cal Z}$.

\begin{cor}\label{3.4}
For every $\rho\in\mathrm{Irr}\,{\cal T}$
one has that $\lambda:=\rho\restriction {\cal Z}\in\mathrm{aut}\,{\cal Z}$.
Then, according to Gelfand's theorem, there exist corresponding
homeomorphisms of $\mathrm{spec}\,{\cal Z}$,
denoted by $f_{\lambda}\in C(\mathrm{spec}\,{\cal Z})$ which are
given by $\lambda(Z)(\phi)=Z(f_{\lambda}^{-1}(\phi)),\,Z\in 
{\cal Z},\,\phi\in \mathrm{spec}\,{\cal Z}$.
\end{cor}

\begin{rem}
Corollary~\ref{3.4} shows that for an irreducible $\rho$
the second possibility considered in \cite[Remark 5.5]{Lledo97b}
of a proper inclusion ${\cal Z}\subset \rho({\cal Z})$
is actually not realized.
\end{rem}

We define the dimension of any object
$\alpha\in\mathrm{Ob}\,{\cal T}$
in the usual way by
$d(\alpha)\EINS:= R_{\alpha}^{\ast}R_{\alpha},$
which satisfies the standard properties of multiplicativity,
additivity etc.~(recall that 
$R_{\alpha}\in (\iota,\overline{\rho}\rho)_{\Mb{C}}$ and
$d(\alpha)>0$). Now using the DR--theory for ${\cal T}_{\Mb{C}}$ 
(cf.~\cite[Sections~2,3]{Doplicher89b}) one arrives at the crucial
decomposition statement for objects.

\begin{pro}\label{Dec}
Let $\alpha\in\mathrm{Ob}\,{\cal T}$. Then
\begin{itemize}
\item[(I)] $d(\alpha)\in \Mb{N},$
\item[(II)] $\alpha= \bigoplus_{j=1}^{r}\bigoplus_{l=1}^{m(\rho_{j},\alpha)}
            \rho_{jl}$, with $\rho_{jl}:=\rho_{j}\in\mathrm{Irr}_{0}\,{\cal T},
            \,l=1,2,\ldots m(\rho_{j},\alpha)$ and \\ 
            $d(\alpha)=\sum_{j=1}^{r}m(\rho_{j},\alpha)d(\rho_{j})$, where
            $m(\rho,\alpha):=\dim\,(\rho,\alpha)_{\Mb{C}}$.
\item[(III)] If $\alpha,\beta\in\mathrm{Ob}\,{\cal T}$, then $\alpha$ 
             is unitarily equivalent to $\beta$ iff 
             $m(\rho,\alpha)=m(\rho,\beta)$ for all $\rho$.
\end{itemize}
\end{pro}

\begin{rem}
Statement (II) in Proposition~\ref{Dec} means explicitly
\[
\alpha(\cdot)=\sum_{\rho,j}W_{\rho j}\,\rho(\cdot)\,W_{\rho j}^{\ast}\,,
\]
where $W_{\rho j}\in(\rho_j,\alpha)_{\Mb{C}}$,
$W_{\rho j}^{\ast}W_{\rho j}=\EINS$, $W_{\rho j}^{\ast}W_{\rho'
j'}=0$, for $(\rho,j)\neq(\rho',j')$ and $\sum_{\rho,j}W_{\rho
j}W_{\rho j}^{\ast}=\EINS$.
Now $\{W_{\rho j}\}_{j}$
is an orthonormal basis of the Hilbert module $(\rho,\alpha)$
and this implies that every orthonormal basis of the Hilbert
module $(\rho,\alpha)$ can be used in the decomposition formula for
$\alpha$. Note that the Hilbert modules $(\rho,\alpha)$
for $\rho\in\mathrm{Irr}_{0}\,{\cal T}$
are mutually orthogonal in ${\cal A}$.
\end{rem}

\subsection{Construction of the Hilbert extension $\{{\cal F},{\cal G}\}$}

In this subsection we will prove part (i) of Theorem~\ref{Main} 
by constructing the Hilbert extension
$\{{\cal F},{\cal G}\}$ of ${\cal A}$ that satisfies
${\cal A}'\cap {\cal F}={\cal Z}$.

We can proceed following the strategy already presented in 
\cite[Sections~3-6]{Baumgaertel97}. To each
$\rho\in\mathrm{Irr}_{0}\,{\cal T}$ we assign a Hilbert space
${\cal H}_{\rho}$ with
$\dim\,{\cal H}_{\rho}=d(\rho)$
and, using orthonormal bases $\{\Phi_{\rho j}\}_{j}$ of
${\cal H}_{\rho}$, we define the ${\cal A}$--left module
\[
{\cal F}_{0}:=\{\sum_{\rho, j}A_{\rho j}\Phi_{\rho j}\mid
               A_{\rho j}\in{\cal A},\,\mbox{finite sum}\},
\]
where the $\{\Phi_{\rho j}\}_{\rho j}$ form an
${\cal A}$--module basis of ${\cal F}_{0}$.
${\cal F}_{0}$ is independent of the special choice of the bases
$\{\Phi_{\rho j}\}_{j}$ of ${\cal H}_{\rho}$ and 
putting $\Phi_{\rho j}A:=\rho(A)\Phi_{\rho j}$,
${\cal F}_{0}$ turns out to be a bimodule.

Further we define Hilbert spaces (recall that 
$\rho < \alpha$ means $\rho$ is a subobject of $\alpha$).
\begin{equation}\label{HS}
{\cal H}_{\alpha}:=\mathop{\oplus}_{\rho<\alpha}(\rho,\alpha)_{\Mb{C}}
\al H._\rho
\quad\mathrm{and}\quad {\cal H}_{\alpha}\subset {\cal
F}_{0},\,\alpha \in\mathrm{Ob}\,{\cal T}\,,
\end{equation}
as well as the right--${\cal Z}$--Hilbert modules
\[
{\ot H.}_{\rho}:={\cal H}_{\rho}{\cal Z}={\cal Z}{\cal H}_{\rho}
\quad\mbox{and}\quad
{\ot H.}_{\alpha}:=\mathop{\oplus}_{\rho<\alpha}
                   (\rho,\alpha){\ot H.}_{\rho}
                  ={\cal H}_{\alpha}{\cal Z},
\]
with the corresponding ${\cal Z}$--scalar product
\[
\langle X,Y\rangle_{\alpha}:= \sum_{\rho,j}\rho^{-1}(X_{\rho j}^{\ast}
Y_{\rho j})\,,\;\mathrm{where}
\]
\[
X=\sum_{\rho,j}X_{\rho j}\Phi_{\rho
j},\,X_{\rho j}\in(\rho,\alpha),\,Y=\sum_{\rho,j}Y_{\rho
j}\Phi_{\rho j},\,Y_{\rho j}\in(\rho,\alpha).
\]

The preceding comments show that we have established the
following functor $\ot F.$ between the categories
${\cal T}$ (resp.~${\cal T}_{\Mb{C}}$)
and the corresponding category of Hilbert ${\cal Z}$--modules
(resp.~Hilbert spaces); (cf.~e.g.~\cite[Section~4]{Lledo97b} and 
\cite[Corollary~3.3]{Baumgaertel97}).

\begin{lem}\label{Funct}
The functor $\ot F.$ given by
\[
\mathrm{Ob}\,{\cal T}\ni\alpha\mapsto {\ot H.}_{\alpha}\subset
{\cal F}_{0}\quad\mathrm{and}\quad(\alpha,\beta)\ni A\mapsto
\ot F.(A)\in{\cal L}_{\cal Z}({\ot H.}_{\alpha}\rightarrow {\ot H.}_{\beta})\,,
\]
where $\ot F.(A)X:=AX,\,X\in{\ot H.}_{\alpha}$,
defines an  isomorphism between the corresponding categories and
$\ot F.(A^{\ast})$ is the module adjoint w.r.t.~$\langle\cdot,\cdot\rangle_{\alpha}$. 
Similarly, one can apply $\ot F.$ to ${\cal T}_{\Mb{C}}$
in order to obtain the associated subcategory of algebraic
Hilbert spaces ${\cal H}_{\alpha}$ and arrows
$\ot F.((\alpha,\beta)_{\Mb{C}})\subset {\cal L}({\cal H}_{\alpha}
\rightarrow {\cal H}_{\beta})$.
\end{lem}
\begin{beweis}
Similar as in \cite[p.~791~ff]{Lledo97b}.
\end{beweis}

Now we can apply the  results in \cite{Baumgaertel97} to the subcategory
$\ot F.({\cal T}_{\Mb{C}})$, in order to enrich gradually the structure  
of ${\cal F}_{0}$:

\begin{lem}\label{Prod}
There exists a product structure on
${\cal F}_{0}$ with the properties
\[
\mathrm{span}\,\{\Phi\cdot\Psi\mid\Phi\in{\cal H}_{\alpha},\,\Psi\in
{\cal H}_{\beta}\}={\cal H}_{\alpha\beta},
\]
\[
\epsilon(\alpha,\beta)\Phi\Psi=\Psi\Phi,\quad \Phi\in{\cal H}_{\alpha},\,
\Psi\in{\cal H}_{\beta},
\]
\[
\langle XY,X'Y'\rangle_{\alpha\beta}=\langle X,X'\rangle_{\alpha}\cdot
\langle Y,Y'\rangle_{\beta},\quad X,X'\in{\cal
H}_{\alpha},\,Y,Y'\in{\cal H}_{\beta}.
\]
\end{lem}

Note that for orthonormal bases
$\{\Phi_{j}\}_{j},\,\{\Psi_{k}\}_{k}$ of
${\cal H}_{\alpha},\,{\cal H}_{\beta}$,
respectively, we obtain from Lemma~\ref{Prod} that
\[
\epsilon(\alpha,\beta)=\sum_{j,k}\Psi_{k}\Phi_{j}
\Psi_{k}^{\ast}\Phi_{j}^{\ast}\,.
\]

As in \cite[Section 5]{Baumgaertel97} we introduce the notion of
a conjugated basis $\Phi_{\overline{\alpha}j}$
of ${\cal H}_{\overline{\alpha}}$ w.r.t.~an orthonormal basis
$\Phi_{\alpha j}$ of ${\cal H}_{\alpha}$ such that
$R_{\alpha}=\sum_{j}\Phi_{\overline{\alpha}j}\Phi_{\alpha j}$.
This is necessary in order to put a compatible *--structure on 
${\cal F}_{0}$.

\begin{lem}
Let $\Phi_{\overline{\rho},j}$ be a conjugated basis corresponding 
to the basis $\Phi_{\rho,j},\,\rho\in\mathrm{Irr}_{0}\,{\cal T}$, and
define $\Phi_{\rho j}^{\ast}:=R_{\rho}^{\ast}\Phi_{\overline{\rho}j},
\,j=1,2,...,d(\rho)$. Then ${\cal F}_{0}$
turns into a *--algebra. The Hilbert spaces ${\cal H}_{\alpha}$
and the corresponding modules ${\ot H.}_{\alpha}$
are algebraic,  i.e.
\[
\langle X,Y\rangle_{\alpha}= X^{\ast}Y,\quad X,Y\in{\ot H.}_{\alpha}.
\]
The objects $\alpha\in\mathrm{Ob}\,{\cal T}$
are identified as canonical endomorphisms
\[
\alpha(A)=\sum_{j=1}^{d(\alpha)}\Phi_{\alpha j}A\Phi_{\alpha j}^{\ast}.
\]
\end{lem}

In ${\cal F}_{0}$ one has natural projections
$\Pi_{\rho}$ onto the $\rho$--component of the decomposition:
\[
\Pi_{\rho}(\sum_{\sigma,j}A_{\sigma j}\Phi_{\sigma j}):=
\sum_{j=1}^{d(\rho)}A_{\rho j}\Phi_{\rho j},\quad
\rho\in\mathrm{Irr}_0\,{\cal T}.
\]
To put a C*--norm
$\Vert\cdot\Vert_{\ast}$
we argue as in \cite[Section~6]{Baumgaertel97}. Its construction 
is essentially based on the following ${\cal A}$--scalar product on
${\cal F}_{0}$
\[
\langle F,G\rangle :=\sum_{\rho,j}\frac{1}{d(\rho)}A_{\rho
j}B_{\rho j}^{\ast},\,\mathrm{for}\; F:=\sum_{\rho,j}A_{\rho j}\Phi_{\rho j},
\,G:=\sum_{\rho,j}B_{\rho j}\Phi_{\rho j},
\]

\begin{lem}
The scalar product $\langle\cdot,\cdot\rangle$ satisfies
$\langle F,G\rangle =\Pi_{\iota}FG^{\ast}$ and
$\Pi_{\rho}$ is selfadjoint w.r.t.~$\langle\cdot,\cdot\rangle$.
The projections $\Pi_{\rho}$
and the scalar product have continuous extensions to
${\cal F}:=\mathrm{clo}_{\Vert\cdot\Vert_{\ast}}{\cal F}_{0}$ and
$\Pi_{\rho}{\cal F}=\mathrm{span}\,\{{\cal A}{\cal H}_{\rho}\}.$
\end{lem}

Finally, the symmetry group w.r.t.~$\langle\cdot,\cdot\rangle$
is defined by the subgroup of all automorphisms
$g\in\mathrm{aut}\,{\cal F}$ satisfying
$\langle gF_{1},gF_{2}\rangle=\langle F_{1},F_{2}\rangle$.
It leads to

\begin{lem} 
The symmetry group coincides with the stabilizer
$\mathrm{stab}\,{\cal A}$ of ${\cal A}$ and the modules
${\ot H.}_{\alpha}$ are invariant w.r.t.~$\mathrm{stab}\,{\cal A}$.
\end{lem}
\begin{beweis}
Use \cite[Lemma~7.1]{Lledo97b}
(cf.~also with the arguments given in \cite[Section~6]{Baumgaertel97}).
\end{beweis}

This suggests to consider the subgroup
${\cal G}\subseteq\mathrm{stab}\,{\cal A}$
consisting of all elements of $\mathrm{stab}\,{\cal A}$
leaving even the Hilbert spaces ${\cal H}_{\alpha}$
invariant. Then it turns out that the pair
$\{{\cal F},{\cal G}\}$
satisfies the properties needed to prove Theorem~\ref{Main}, 
i.e.~$\{{\cal F},{\cal G}\}$ is a Hilbert extension of
${\cal A}$. The following result concludes the proof of
part~(i) of Theorem~\ref{Main}.

\begin{lem}
${\cal G}$ is compact and the spectrum $\mathrm{spec}\,{\cal G}$ on
${\cal F}$ coincides with the dual $\hat{\cal G}$. For
$\rho\in\mathrm{Irr}\,{\cal T}$ the Hilbert spaces
${\cal H}_{\rho}$ are irreducible w.r.t.~${\cal G},$
i.e.~there is a bijection
$\mathrm{Irr}_0\,{\cal T}\ni\rho\leftrightarrow D\in\hat{\cal G}$.
Moreover ${\cal A}$ coincides with the fixed point algebra 
of the action of ${\cal G}$ in ${\cal F}$ and
${\cal A}'\cap {\cal F}={\cal Z}$.
\end{lem}

\begin{rem}
\begin{itemize}
\item[(i)] From \cite[Section~7]{Lledo97b} it follows that 
 $\mathrm{stab}\,{\cal A}$ is in general {\em not} compact.
\item[(ii)] The characterization of $\mathrm{stab}\,{\cal A}$
given in \cite[Theorem~7.11]{Lledo97b} in terms of functions contained
in $C(\mathrm{spec}\al Z.\rightarrow\al G.)$ is in general not correct, 
although in some special cases like the one--dimensional torus 
${\cal G}:=\Mb{T}$ it is true that
$\mathrm{stab}\,{\cal A}\cong C(\mathrm{spec}\al Z.\rightarrow \Mb{T})$
(cf.~\cite{pBaumgaertel00}). It is though possible to give a similar
characterization of $\mathrm{stab}\al A.$ in terms of functions
contained in $C(\mathrm{spec}\al Z.\rightarrow\mathrm{Mat}(\C))$.
\end{itemize}
\end{rem}

\subsection{Uniqueness result}

Now we prove part~(ii) of Theorem~\ref{Main}. 
First assume that the subcategories
${\cal T}_{\Mb{C}}$ and ${\cal T}_{\Mb{C}}'$
are equivalent. We consider the Hilbert extension ${\cal F}$
assigned to ${\cal T}_{\Mb{C}}$.
The corresponding invariant Hilbert spaces are given by (\ref{HS}). 
Now we change these Hilbert spaces by
\[
{\cal H}_{\alpha}\rightarrow V_{\alpha}{\cal H}_{\alpha}=:\al H._\alpha'\,.
\]
Using the function $\ot F.$ of Lemma~\ref{Funct} so that
${\cal L}_{\cal G}({\cal H}_{\alpha}\rightarrow{\cal H}_{\beta})
:=\ot F.((\alpha,\beta)_{\Mb{C}})\cong (\alpha,\beta)_{\Mb{C}}$
we obtain
\begin{equation}\label{Intw}
{\cal L}_{\cal G}(V_{\alpha}{\cal H}_{\alpha}\rightarrow
V_{\beta}{\cal H}_{\beta})=
V_{\beta}{\cal L}_{\cal G}({\cal H}_{\alpha}\rightarrow
{\cal H}_{\beta})V_{\alpha}^{\ast}\cong
(\alpha,\beta)_{\Mb{C}}'.
\end{equation}
Further, w.r.t.~the ``new Hilbert spaces" we obtain the `primed'
permutators and conjugates of the second subcategory. This means, it
is sufficient to prove that if the subcategory
${\cal T}_{\Mb{C}}$
is given, then two Hilbert extensions, assigned to
$({\cal T},{\cal T}_\C)$
according to the first part of the theorem, are always
${\cal A}$--module isomorphic. Now let
${\cal F}_{1},{\cal F}_{2}$
be two Hilbert extensions assigned to
${\cal T}_{\Mb{C}}$.
For $\rho\in\mathrm{Irr}_{0}\,{\cal T}$ let
$\{\Phi_{\rho j}^{1}\}_j$,$\{\Phi_{\rho j}^{2}\}_k$
be orthonormal bases of the Hilbert spaces
${\cal H}_{\rho}^{1},{\cal H}_{\rho}^{2}$,
respectively. Then
\[
\Phi_{\rho j}^{r}\cdot\Phi_{\sigma k}^{r}=
\sum_{\tau, l}K_{\rho j\sigma k}^{\tau l}\Phi_{\tau l}^{r},\quad
K_{\rho j\sigma k}^{\tau l}\in(\tau,\rho\sigma)_{\Mb{C}},\,r=1,2.
\]
Therefore the definition
\[
\al J.(\sum_{\rho,j}A_{\rho j}\Phi_{\rho j}^{1}):=
\sum_{\rho,j}A_{\rho j}\Phi_{\rho j}^{2}
\]
is easily seen to extend to an
${\cal A}$--module isomorphism from
${\cal F}_{1}$ onto ${\cal F}_{2}$ (see 
\cite[p.~203~ff.]{bBaumgaertel92}).

Second, we assume that the Hilbert extensions
${\cal F}_{1},{\cal F}_{2}$ assigned to
${\cal T}_{\Mb{C}}^{1},{\cal T}_{\Mb{C}}^{2}$,
respectively, are 
${\cal A}$--module isomorphic. The 
${\cal G}$--invariant Hilbert spaces are given by (\ref{HS}). 
Now let $\al J.$ be an ${\cal A}$--module isomorphism 
$\al J.\colon\ {\cal F}_{1}\rightarrow {\cal F}_{2}$ so that
\[
\al J.({\cal H}_{\alpha}^{1})=\mathop{\oplus}_{\rho<\alpha}
(\rho,\alpha)_{\Mb{C}}^{1}\,\al J.({\cal H}_{\rho}^{1})
\]
and again the $\al J.({\cal H}_{\alpha}^{1})$
form a system of ${\cal G}$--invariant Hilbert spaces in 
${\cal F}_{2}$. Further we have the system
${\cal H}_{\alpha}^{2}$ in ${\cal F}_{2}$.
That is, to each $\alpha$ we obtain two
${\cal G}$--invariant Hilbert spaces ${\cal H}_{\alpha}^{2}$
and $\al J.({\cal H}_{\alpha}^{1})$
that are contained in the Hilbert module $\ot H._\alpha^2$.
Let $\{\Phi_{\alpha,j}\}_j$,$\{\Psi_{\alpha,j}\}_j$ be orthonormal
bases of $\al J.({\cal H}_{\alpha}^{1})$,$\al H._\alpha^2$,
respectively. Then obviously 
$V_\alpha:=\sum_j \Psi_{\alpha,j}\Phi_{\alpha,j}^*$ is a unitary
with $V_\alpha\in(\alpha,\alpha)$ and $\al H._\alpha^2=
V_\alpha\,\al J.({\cal H}_{\alpha}^{1})$. Further, 
for $X\in\al H.^1_\alpha$,$Y\in\al H.^1_\beta$ (hence 
$XY\in\al H.^1_{\alpha\beta}$) we have 
\[
 V_\alpha \al J.(X) V_\beta \al J.(Y)
    = V_\alpha\alpha(V_\beta)\,\al J.(XY) 
    = V_{\alpha\circ\beta}\,\al J.(XY)\,,
\]
and this implies $V_{\alpha\circ\beta}=V_{\alpha}\times V_{\beta}$.
Finally, we argue as in (\ref{Intw}) to obtain
\[
 V_\beta\, (\alpha,\beta)_\C^1 \,V_\alpha=(\alpha,\beta)_\C^2.
\]
and the latter equation implies Eqs.~(\ref{Int})--(\ref{RS}).

\section{Conclusions}
In the present paper we present the solution of the problem of 
finding the unique (up to $\al A.$--moldule isomorphy) 
Hilbert extension $\{\al F.,\al G.\}$ 
of a unital C*--algebra $\al A.$ with nontrivial center
$\al Z.$, given a suitable endomorphism category $\al T.$ of 
$\al A.$, which we characterize in Section~\ref{Assu}.
The extension satisfies $\al A.'\cap\al F.=\al Z.$ and the essential
step for its construction is the specification of a subcategory
$\al T._\C$ of $\al T.$, which is of the well--known DR--type.
From the point of view of the DR--theory the appearence of the 
subcategory $\al T._\C\subset\al T.$ is quite natural, since the group 
appearing in the extension is still {\em compact}. There are several
directions in which the present results could be generalized. First,
we hope that the inclusion situation $\al T._\C\subset\al T.$
may also be relevant for braided tensor categories, since in this
context there are 2--dimensional physically relevant models where
a nontrivial center appears (see e.g.~\cite{Buchholz88,Fredenhagen92}).
Second, one could try to find extensions, where the condition on the 
relative commutant $\al A.'\cap\al F.=\al Z.$ (which is crucial 
for our approach) is not satisfied anymore. In this context 
non unitarily equivalent
irreducible endomorphisms will no longer be disjoint and one needs
probably to replace the free modules $\ot H.$ that appeared in our
approach by more general C*--Hilbert modules. Finally, we hope that the 
present results as well as those in \cite{Lledo97b} will motivate
a more systematic study of the representation theory of (say 
compact) groups over Hilbert C*--modules.


\providecommand{\bysame}{\leavevmode\hbox to3em{\hrulefill}\thinspace}

\end{document}